\documentclass[12pt]{article}
\usepackage{latexsym,amsfonts,amssymb,mathrsfs}
\usepackage[dvips]{color}
\usepackage{colordvi,multicol}

\newtheorem{thm}{Theorem}[section]
\newtheorem{cor}[thm]{Corollary}
\newtheorem{lem}[thm]{Lemma}

\newtheorem{defi}[thm]{Definition}
\newtheorem{rem}[thm]{Remark}

\begin{document}
\title{\bf Unlinking Theorem for Symmetric
Quasi-convex Polynomials}
\author{He-Jing Hong$^{a,b}$, Ze-Chun Hu$^{c,}$\thanks{{Corresponding
author: College of Mathematics, Sichuan University, Chengdu 610065, China.} \vskip 0cm {\it \ E-mail address:} zchu@scu.edu.cn}\\ \\
{\small $^a$Department of Mathematics, Nanjing University}\\
{\small $^b$Clinchoice Inc. of Nanjing}\\
{\small $^c$College of Mathematics, Sichuan University}}

\maketitle
\date{}

% \centerline{Ze-Chun Hu\thanks{Corresponding author.}} \centerline{\small College
%of Mathematics, Sichuan University, Chengdu, 610064, China}
%\centerline{\small E-mail: zchu@scu.edu.cn}
%\vskip 0.7cm \centerline{Qian-Qian Zhou} \centerline{\small Department of
%Mathematics, Nanjing University, Nanjing, 210093, China} \centerline{\small
%E-mail: 15266479708@163.com}

%\vskip 1.4cm

\vskip 0.5cm \noindent{\bf Abstract}\quad Let $\mu_n$ be the standard Gaussian
measure on $\mathbb{R}^n$ and $X$ be a random vector on
$\mathbb{R}^n$ with the law $\mu_n$. U-conjecture states that if $f$
and $g$ are two polynomials on $\mathbb{R}^n$ such that $f(X)$ and
$g(X)$ are independent, then there exist an orthogonal
transformation $L$ on $\mathbb{R}^n$ and an integer $k$ such that
$f\circ L$ is a function of $(x_1,\cdots,x_k)$ and $g\circ L$ is a
function of $(x_{k+1},\cdots,x_n)$. In this case, $f$ and $g$ are
said to be unlinked. In this note, we prove that two symmetric,
quasi-convex polynomials $f$ and $g$ are unlinked if $f(X)$ and
$g(X)$ are independent.\\

\smallskip

\noindent {\bf Keywords}\quad U-conjecture; quasi-convex polynomial;
Gaussian correlation conjecture

\smallskip

\noindent {\bf Mathematics Subject Classification (2010)}\quad 60E15; 62H05

%\tableofcontents

\section{Introduction and main result}

Let $\mu_n$ be the standard Gaussian measure on $\mathbb{R}^n(n\geq 2)$ and
$X$ be a random vector on $\mathbb{R}^n$ with the law $\mu_n$. In
1973, Kagan, Linnik and Rao \cite{KLR73} considered the following
problem: if $f$ and $g$ are two polynomials on $\mathbb{R}^n$ such
that $f(X)$ and $g(X)$ are independent, then is it possible to find an
orthogonal transformation $L$ on $\mathbb{R}^n$ and an integer $k$
such that $f\circ L$ is a function of $(x_1,\cdots,x_k)$ and $g\circ
L$ is a function of $(x_{k+1},\cdots,x_n)$? If the answer is
positive, then
 $f$ and $g$ are said to be unlinked. This problem is called U-conjecture and is still open.

The U-conjecture is true for the case $n=2$, and some special cases
have been proved for larger number of variables (see Sections
11.4-11.6 of \cite{KLR73}). In 1994, Bhandari and DasGupta
\cite{BD94} proved that the U-conjecture holds for two symmetric convex functions $f$ and $g$  under an additional condition. The additional condition can be canceled since
the Gaussian correlation conjecture has been proved (see Royen \cite{Royen} or Lata{\l}a and Matlak \cite{LM}).

 Bhandari and Basu \cite{BB06} proved that the U-conjecture holds for two nonnegative convex polynomials $f$ and $g$ with $f(0)=0$.  Harg\'{e} \cite{Har05} proved that if $f,g: \mathbb{R}^n\to \mathbb{R}$ are two convex functions in
$L^2(\mu_n)$, and $f$ is a real analytic function satisfying
$f(x)\geq f(0),\forall x\in \mathbb{R}^n$, and $f$ and $g$ are
independent with respect to $\mu_n$, then they are unlinked.

Malicet et al. \cite{MNPP16} proved that the U-conjecture is true when $f,g$ belong to a class of polynomials, which is defined based on the infinitesimal generator of Ornstein-Uhlenbeck semigroup.

In Remark 2 of \cite{BB06}, the authors wish that their  result
could be extended to symmetric, quasi-convex polynomials. In this
note, we will give an affirmative answer based on the first author's master thesis \cite{Hong09} and  prove  the following result.

\begin{thm}\label{main-thm}
Two symmetric, quasi-convex polynomials $f$ and $g$ are unlinked if $f$ and $g$ are independent with respect to $\mu_n$.
\end{thm}

\section{Proof of Theorem \ref{main-thm}}\setcounter{equation}{0}

Before giving the proof of Theorem \ref{main-thm}, we present some preliminaries.

A function $f:\mathbb{R}^n\to \mathbb{R}$ is called {\it
quasi-convex} if for any $\alpha\in [0,1]$ and any $x,y\in
\mathbb{R}^n$,
$$
f(\alpha x+(1-\alpha)y)\leq \max\{f(x),f(y)\}.
$$
It's easy to know that a convex function is quasi-convex. About the
properties of quasi-convex functions, and the relations between
convex and quasi-convex functions,  refer to a survey paper
Greenberg and Pierskalla \cite{GP71}.

\begin{lem}\label{lem-2.1}
Suppose that $g:\mathbb{R}\to \mathbb{R}$ is a quasi-convex
polynomial and there exist $\lambda_1,\lambda_2\in \mathbb{R}$ such
that $g(\lambda_1)\neq g(\lambda_2)$. Then one of the following two
claims holds.\\
(a) There exists $\lambda_0$ such that $g(u)<g(v)$ for any
$\lambda_0\leq u<v$
 and
 $\lim\limits_{\lambda\to\infty}g(\lambda)=\infty.$\\
(b) There exists $\lambda_0$ such that $g(u)<g(v)$ for any $v<u\leq
\lambda_0$ and
 $\lim\limits_{\lambda\to -\infty}g(\lambda)=\infty.$
\end{lem}
{\bf Proof.} Since $g$ is a polynomial on $\mathbb{R}$, we can write
it as
\begin{eqnarray}\label{lem-1-a}
g(\lambda)=a_n\lambda^n+a_{n-1}\lambda^{n-1}+\cdots+ a_1\lambda+a_0.
\end{eqnarray}
By the assumption, $g$ is not a constant, so $n\geq 1$ and $a_n\neq
0$. By (\ref{lem-1-a}), we obtain
\begin{eqnarray}\label{lem-1-b}
g'(\lambda)=na_n\lambda^{n-1}+(n-1)a_{n-1}\lambda^{n-2}+\cdots +a_1.
\end{eqnarray}

Without loss of generality, we can  assume that
$\lambda_1<\lambda_2$. We have the following two cases:

{\bf Case 1:} $g(\lambda_1)<g(\lambda_2)$. Define
$h(\lambda):=g(\lambda)- g(\lambda_1)$. Then $h(\lambda_1)=0$. By the
definition of quasi-convex function, we  know that $h$ is
quasi-convex, and for any $\lambda>\lambda_2$, we have
\begin{eqnarray}\label{lem-1-b-1}
h(\lambda_2)&=&h\left(\frac{\lambda-\lambda_2}{\lambda-\lambda_1}\lambda_1+
\frac{\lambda_2-\lambda_1}{\lambda-\lambda_1}\lambda\right)\leq\max\{h(\lambda_1),h(\lambda)\}=\max\{0,h(\lambda)\}.
\end{eqnarray}
Since  $h(\lambda_2)=g(\lambda_2)-g(\lambda_1)>0$, by
(\ref{lem-1-b-1}), we ge that for any $\lambda>\lambda_2$,
$h(\lambda_2)\leq h(\lambda)$, i.e.
\begin{eqnarray}\label{lem-1-c}
g(\lambda_2)\leq g(\lambda),\ \forall \lambda>\lambda_2.
\end{eqnarray}
By (\ref{lem-1-a}) and (\ref{lem-1-c}), we get that $a_n>0$, and
thus
\begin{eqnarray}\label{lem-1-c-1}
&&\lim_{\lambda\to \infty}g(\lambda)=\lim_{\lambda\to \infty}
\left(a_n\lambda^n+a_{n-1}\lambda^{n-1}+\cdots+
a_1\lambda+a_0\right)=\infty.
\end{eqnarray}

If $n=1$, then $g(\lambda)=a_1\lambda+a_0$ with $a_1>0$, and thus (a) holds in this case.
If $n\geq 2$, then
\begin{eqnarray}\label{lem-1-d}
&& \lim_{\lambda\to \infty}g'(\lambda)=
 \lim_{\lambda\to \infty}\left(na_n\lambda^{n-1}+(n-1)a_{n-1}\lambda^{n-2}+\cdots +a_1\right)=\infty.
  \end{eqnarray}
By (\ref{lem-1-d}), there exists $\lambda_0$ such that for any
$\lambda>\lambda_0$, $g'(\lambda)>0$, which together with  (\ref{lem-1-c-1}) implies that (a) holds in this case.

{\bf Case 2:} $g(\lambda_1)>g(\lambda_2)$. Define
$\bar{h}(\lambda):=g(\lambda)- g(\lambda_2)$. Then
$\bar{h}(\lambda_2)=0$, and as in Case 1, $\bar{h}$ is a
quasi-convex function and for any $\lambda<\lambda_1$, we have
\begin{eqnarray}\label{lem-1-d-1}
\bar{h}(\lambda_1)&=&\bar{h}\left(\frac{\lambda_1-\lambda}{\lambda_2-\lambda}\lambda_2+
\frac{\lambda_2-\lambda_1}{\lambda_2-\lambda}\lambda\right)
\leq\max\{\bar{h}(\lambda_2),\bar{h}(\lambda)\}=\max\{0,\bar{h}(\lambda)\}.
\end{eqnarray}
Since $\bar{h}(\lambda_1)=g(\lambda_1)-g(\lambda_2)>0$, by
(\ref{lem-1-d-1}), we obtain that for any $\lambda<\lambda_1$,
$\bar{h}(\lambda_1)\leq \bar{h}(\lambda)$, i.e.
\begin{eqnarray}\label{lem-1-e}
g(\lambda_1)\leq g(\lambda),\ \forall \lambda<\lambda_1.
\end{eqnarray}
By (\ref{lem-1-e}) and (\ref{lem-1-a}), we know that one of the
following two claims must hold:

 (i) $n$ is  even and  $a_n>0$;\quad\quad\quad (ii) $n$ is  odd  and  $a_n<0$.
\bigskip

If (i) holds, then by the proof of Case 1 above, we know that (a) is
true.

If (ii) holds, then
\begin{eqnarray}\label{lem-1-e-1}
&&\lim_{\lambda\to -\infty}g(\lambda)=\lim_{\lambda\to -\infty}
\left(a_n\lambda^n+a_{n-1}\lambda^{n-1}+\cdots+
a_1\lambda+a_0\right)=\infty.
\end{eqnarray}
If $n=1$, then $g(\lambda)=a_1\lambda+a_0$ with $a_1<0$, and thus (b) holds in this case.
If $n\geq 3$, then
\begin{eqnarray}\label{lem-1-f}
&& \lim_{\lambda\to -\infty}g'(\lambda)=
 \lim_{\lambda\to -\infty}\left(na_n\lambda^{n-1}+(n-1)a_{n-1}\lambda^{n-2}+\cdots +a_1\right)=-\infty.
  \end{eqnarray}
By (\ref{lem-1-f}), there exists $\lambda_0$ such that for any
$\lambda<\lambda_0$, $g'(\lambda)<0$, which together with (\ref{lem-1-e-1}) implies that  (b) holds in this case.
\hfill\fbox

\begin{cor}\label{cor-2.2}
Let $g:\mathbb{R}\to \mathbb{R}$ be a quasi-convex polynomial. If
$g$ has an upper bound, then $g$ is a constant function.
\end{cor}

\begin{cor}\label{cor-2.3}
Let $U: \mathbb{R}^n\to\mathbb{R}$ be a quasi-convex polynomial.
Suppose that for  two fixed vectors
$\beta_1,\beta_2\in\mathbb{R}^n$, $U\left(
\beta_1+\lambda\beta_2\right)$ is a constant funciton of $\lambda\in
\mathbb{R}$. Then for any fixed vector $b\in\mathbb{R}^n$,
$U\left(b+\lambda\beta_2 \right)$ is a constant function of
$\lambda$.
\end{cor}

\noindent{\bf Proof. }  For any fixed vector $b\in \mathbb{R}^n$,
define $g(\lambda)=U(b+\lambda \beta_2),\lambda\in \mathbb{R}$. Then
$g(\lambda)$ is a polynomial of $\lambda$. By the quasi-convexity of
$U$, we know that for any $\alpha\in [0,1]$ and
$\lambda_1,\lambda_2\in \mathbb{R}$, we have
\begin{eqnarray*}
g(\alpha \lambda_1+(1-\alpha)\lambda_2)&=&U(b+\left(\alpha
\lambda_1+(1-\alpha)\lambda_2\right) \beta_2)\\
&=&U(\alpha(b+\lambda_1\beta_2)+(1-\alpha)(b+\lambda_2\beta_2))\\
&\leq&\max\{U(b+\lambda_1\beta_2),U(b+\lambda_2\beta_2)\}\\
&=&\max\{g(\lambda_1),g(\lambda_2)\}.
\end{eqnarray*}
Thus $g(\lambda)$ is a quasi-convex polynomial. By the
quasi-convexity of $U$,
\begin{eqnarray}\label{cor2.3-1}
 g(\lambda)=U(
b+\lambda\beta_2)&=&U\left(\frac{1}{2}(2b-\beta_1)+\frac{1}{2}(\beta_1+2\lambda
\beta_2)\right)\nonumber\\
&\leq& \max\left\{U\left( 2b-\beta_1\right),U\left(
\beta_1+2\lambda\beta_2\right)\right\}.
\end{eqnarray}
By (\ref{cor2.3-1}) and the assumption that $U\left(
\beta_1+\lambda\beta_2\right)$ is a constant funciton of $\lambda$,
we get that the quasi-convex polynomial $g(\lambda)$ has an upper
bound. Hence by Corollary \ref{cor-2.2}, we know that
$U\left(b+\lambda\beta_2 \right)$ is a constant function of
$\lambda$. \hfill\fbox

\begin{cor}\label{lem-2.4}
Let $U: \mathbb{R}^n\to\mathbb{R}$ be a quasi-convex polynomial with
$U\left(0 \right)=0$. Define
\begin{eqnarray}\label{lem-2.4-a}
S_U:=\left\{\alpha:U\left(\lambda\alpha
\right)=0,\ \forall \lambda\in\mathbb{R} \right\}.
\end{eqnarray}
Then $S_U$ is a vector subspace of $\mathbb{R}^n$.
\end{cor}
\noindent {\bf  Proof.} Let $\alpha_1,\alpha_2\in S_U$. For any $
c_1,c_2,\lambda\in\mathbb{R}$, by Corollary \ref{cor-2.3} we get
that
\begin{eqnarray*}
U\left(\lambda(c_1\alpha_1+c_2\alpha_2) \right)=U\left(\lambda
c_1\alpha_1+\lambda c_2\alpha_2 \right)=U\left(\lambda c_1\alpha_1
\right)=0.
\end{eqnarray*}
Hence $c_1\alpha_1+c_2\alpha_2 \in S_U$, and thus $S_U$ is a vector subspace of $\mathbb{R}^n$. \hfill\fbox

\bigskip

Now suppose that $U$ and $V$ are two quasi-convex polynomials from
$\mathbb{R}^n$ into $\mathbb{R}$ satisfying that $U(0)=V(0)=0$.
Define $S_U$ by (\ref{lem-2.4-a}). Similarly, define  $S_V$.

\begin{defi}
$U$ and $V$ are said to be concordant of order $r$, if
\begin{eqnarray}\label{def}
\dim(S_U^{\bot})-\dim(S_U^{\bot}\cap S_V)=r.
\end{eqnarray}
\end{defi}

Note that this definition is symmetric in $U$ and $V$, i.e. if
(\ref{def}) holds, then (see \cite{BD94})
\begin{eqnarray*}
\dim(S_V^{\bot})-\dim(S_V^{\bot}\cap S_U)=r.
\end{eqnarray*}

\begin{thm}\label{thm}
 Let $X$ be an $n\times 1$ random vector distributed as $N\left(0,I_n
\right)$. Let $U$ and $V$ be two symmetric $(i.e.\ U(x)=U(-x),
V(x)=V(-x))$ quasi-convex polynomials on $\mathbb{R}^n$ satisfying
$Cov\left(U\left(X \right),V\left(X\right) \right)=0$. Furthermore,
assume that $U\left(0 \right)=V\left(0 \right)=0$, and $U$ and $V$
are concordant of order $r$. Then there exists an orthogonal
transformation $Y=LX$ such that $U$ and $V$ can be expressed as
functions of two different sets of components of $Y$, i.e. $U$ and
$V$ are unlinked.
\end{thm}

\noindent {\bf Proof.} Based on the lemmas and corollaries
established above, the proof of this theorem is similar to the one
of \cite{BD94}. For the reader's convenience, we spell out the details in the following.

Let $\left\{\alpha_1,...,\alpha_{r+t} \right\}$,
$\left\{\alpha_{r+1},...,\alpha_{r+t}  \right\}$,
$\left\{\alpha_1,...,\alpha_{r+t+m} \right\}$ and
$\left\{\alpha_1,...,\alpha_{n} \right\}$ be orthonormal bases of
$S_U^\bot$,  $S_U^\bot\bigcap S_V$ , $S_U^\bot+S_V^\bot$, and
$\mathbb{R}^n$, respectively. We will show that if $r>0$ then
$Cov(U(X),V(X))>0$, which contradicts the condition given in the
theorem, and so we get $r=0$, and thus $U$ and $V$ are unlinked.

Define $Y_1,Y_2,\cdots,Y_n$ by $X=\sum_{i=1}^{n}{Y_i\alpha_i}$, i.e.
$Y_i$ is the $i$-th component of $X$. Then $Y_1,Y_2,\cdots,Y_n$ are
i.i.d. as $N\left(0,1\right)$. By Corollary \ref{cor-2.3},
\begin{eqnarray*}
&&U\left(X \right)=U\left(\sum_{i=1}^{n}{Y_i\alpha_i}
\right)=U\left(\sum_{i=1}^{r}{Y_i\alpha_i}+\sum_{i=r+1}^{r+t}{Y_i\alpha_i}
\right),\\
&&V\left(X \right)=V\left(\sum_{i=1}^{n}{Y_i\alpha_i}
\right)=U\left(\sum_{i=1}^{r}{Y_i\alpha_i}+\sum_{i=r+t+1}^{r+t+m}{Y_i\alpha_i}
\right).
\end{eqnarray*}

Assume  that $r>0$. Let $y^*=\left(y_1,...,y_r \right)'$  be a
nonzero vector in $\mathbb{R}^r$. Define
\begin{eqnarray*}
&&U^*\left(y^*
\right):=E\left[U\left(\sum_{i=1}^{r}{y_i\alpha_i}+\sum_{i=r+1}^{r+t}{Y_i\alpha_i}
\right) \right],\\
&& V^*\left(y^*
\right):=E\left[V\left(\sum_{i=1}^{r}{y_i\alpha_i}+\sum_{i=r+t+1}^{r+t+m}{Y_i\alpha_i}
\right) \right].
\end{eqnarray*}
 Then $U^*$ and $V^*$ are two symmetric quasi-convex polynomials of
 $y^*$.

By the choice of the bases, $U\left(\lambda\sum_{i=1}^{r}{y_i
\alpha_i} \right)$ is not a zero function of $\lambda$. By Corollary
\ref{cor-2.3} and the condition $U(0)=0$, we know that $U(\lambda\sum_{1}^{r}{y_i x_i} +$
$\sum_{r+1}^{r+t}{y_i x_i})$ is not a constant of $\lambda$. In
addition, by the symmetry and quasi-convexity of $U$, $U(x)\geq
U(0),\forall x\in \mathbb{R}^n$. Hence by Lemma \ref{lem-2.1}, we
get that when $\lambda\rightarrow\infty$,
\begin{eqnarray}\label{thm-2.7-0}
U\left(\lambda\sum_{1}^{r}{y_i \alpha_i} +\sum_{r+1}^{r+t}{Y_i
x_i}\right)+U\left(-\lambda\sum_{1}^{r}{y_i \alpha_i}
+\sum_{r+1}^{r+t}{Y_i x_i}\right)\rightarrow\infty.
\end{eqnarray}
Taking the expectation of (\ref{thm-2.7-0}) with respect to
$Y_{i+1},...,Y_{r+t}$ and using Egoroff's theorem (see e.g. \cite[Theorem 21.3]{Mu52} or \cite[Remark 2.3.6(1)]{Yan04}), we obtain
\begin{eqnarray}\label{thm-2.7-a}
U^*\left(\lambda y^* \right )\rightarrow\infty\quad\mbox{as}\quad
\lambda\to\infty.
\end{eqnarray}
Similarly,
\begin{eqnarray}\label{thm-2.7-b}
V^*\left(\lambda y^* \right )\rightarrow\infty\quad \mbox{as}\quad
\lambda\to\infty.
\end{eqnarray}

Define $Y^*=(Y_1,\cdots,Y_r)'$. By the independence of components of
$X=(Y_1,\cdots,Y_r,Y_{r+1},\cdots,Y_n)'$ and simple calculations, we
have
\begin{eqnarray}\label{thm-2.7-c}
&&Cov\left(U\left(X \right),V\left( X\right) \right)\nonumber\\
&&=
E[U\left(X \right)V\left(X \right)]-E[U\left(X \right)]E[V\left(X \right)]\nonumber\\
&&=E\left[U^*\left(Y^* \right)V^*\left(Y^* \right)\right]-
E\left[U^*\left(Y^* \right)\right]E\left[V^*\left(Y^* \right)\right]\nonumber\\
&&=\int_{0}^{\infty}{\int_{0}^{\infty}{\left[P\left(Y^*\in
A_{k_1}^c\cap B_{k_2}^c\right) -
P\left(Y^*\in A_{k_1}^c \right)P\left(Y^*\in B_{k_2}^c \right)\right]dk_1}dk_2}\nonumber\\
&&=\int_{0}^{\infty}{\int_{0}^{\infty}{\left[P\left(Y^*\in
A_{k_1}\cap B_{k_2}\right) - P\left(Y^*\in A_{k_1}
\right)P\left(Y^*\in B_{k_2} \right)\right]dk_1}dk_2},
\end{eqnarray}
where
\begin{eqnarray*}
A_{k_1}=\left\{y^*:U^*\left(y^* \right)\leq k_1 \right\},\ \ \
B_{k_2}=\left\{y^*:V^*\left(y^* \right)\leq k_2 \right\}.
\end{eqnarray*}

Since $U^*(y^*)$ and $V^*(y^*)$ are symmetric, quasi-convex
polynomials of $y^*$, ${A_k}_1$ and $ {B_k}_2$ are both symmetric
convex sets (see  \cite[Table II]{GP71}). By the Gaussian correlation inequality (see \cite{Royen} or \cite{LM}),
\begin{eqnarray}\label{thm-2.7-d}
P\left(Y^*\in A_{k_1}\cap B_{k_2}\right) - P\left(Y^*\in A_{k_1}
\right)P\left(Y^*\in B_{k_2} \right)\geq 0.
\end{eqnarray}

Define a set
\begin{eqnarray*}
M=\{(k_1,k_2)\in (0,\infty)\times (0,\infty)|A_{k_1}\subset
B_{k_2},P(Y^*\in B_{k_2}^c)>0,P(Y^*\in A_{k_1})>0\}.
\end{eqnarray*}
When $A_{k_1}\subset B_{k_2}$, we have
\begin{eqnarray*}
&&P(Y^*\in A_{k_1}\cap B_{k_2}) - P(Y^*\in A_{k_1} )P(Y^*\in
B_{k_2})\\
&&=P(Y^*\in A_{k_1})(1-P(Y^*\in B_{k_2}))\\
&&=P(Y^*\in A_{k_1})P(Y^*\in B_{k_2}^c).
\end{eqnarray*}
Hence we obtain
\begin{eqnarray}\label{thm-2.7-e}
M\subset \left\{(k_1,k_2)\in (0,\infty)\times (0,\infty)|P(Y^*\in
A_{k_1}\cap B_{k_2}) - P(Y^*\in A_{k_1} )P(Y^*\in
B_{k_2})>0\right\}.
\end{eqnarray}
By (\ref{thm-2.7-a}), (\ref{thm-2.7-b}), and Lemma \ref{lem-2.1},
the Lebesgure measure of $M$ is positive. Hence by (\ref{thm-2.7-c})
, (\ref{thm-2.7-d}) and (\ref{thm-2.7-e}), we obtain
$$
Cov(U(X),V(X))>0,
$$
which contradicts the assumption, and so $r=0$.\hfill\fbox

\begin{rem}
The assumption that $U(0)=V(0)=0$ in Theorem \ref{thm} can be taken
out since by the symmetry and quasi-convexity of $U$ and $V$, we
have $U(x)\geq U(0),V(x)\geq V(0)$ for all $x\in \mathbb{R}^n$, and
we can consider the polynomials $U(x)-U(0)$ and $V(x)-V(0)$, which
satisfy the conditions in Theorem \ref{thm}.
\end{rem}

\vskip 0.5cm
{ \noindent {\bf\large Acknowledgments}

\noindent  This work was supported by National Natural Science Foundation of China (Grant Nos. 11771309; 11871184).

\end{document}